\documentclass[10pt]{amsart}
\usepackage{amssymb}
\usepackage{pstricks}
\usepackage[dvips]{graphicx}

\newtheorem{theorem}{Theorem}

\newtheorem{lemma}{Lemma}
\newtheorem{sublemma}{Sublemma}
\newtheorem{claim}{Claim}
\newtheorem{corollary}{Corollary}

\theoremstyle{definition}
\newtheorem{definition}{Definition}
\newtheorem{example}{Example}

\begin{document}

\title{On the local structure of doubly laced crystals}

\author[P.~Sternberg]{Philip Sternberg}
\address{Department of Mathematics, University of California, One Shields
Avenue, Davis, CA 95616-8633, U.S.A.}
\email{sternberg@math.ucdavis.edu}
\urladdr{http://www.math.ucdavis.edu/$\sim$sternberg}

\thanks{\textit{Date:} February 2006}
\thanks{Supported in part by NSF grants DMS-0135345, DMS-0200774, and 
DMS-0501101.}


\begin{abstract}
    Let $\mathfrak{g}$ be a Lie algebra all of whose regular
    subalgebras of rank 2 are type $A_{1}\times A_{1}$, $A_{2}$, or
    $C_{2}$, and let $B$ be a crystal graph corresponding to a
    representation of $\mathfrak{g}$.  We explicitly describe the
    local structure of $B$, confirming a conjecture of Stembridge.
\end{abstract}

\maketitle

\section{Introduction}

Since their introduction by Kashiwara \cite{Kash:1990, Kash:1991},
crystal bases have proven very useful in the study of representation
theory.  In particular, given any highest weight integrable module $V$
over a symmetrizable quantum group we can construct a colored directed
graph, called a crystal graph, that encodes nearly all the
representation theoretic information of $V$.  Alternatively, one may
define crystals axiomatically; many examples of axiomatic crystals
that do not correspond to any representation of a quantum group are
known.

Many explicit combinatorial models have been developed for crystal
graphs of representations; two examples are paths in the weight space
of the algebra being represented \cite{L1,L2} and generalized Young
tableaux \cite{KN}.  In all such constructions, the combinatorics of
the crystals are defined by global poroperties.  In \cite{Stem},
Stembridge introduced a set of graph theoretic axioms, each of which
addresses only local properties of a colored directed graph, that
characterizes highest weight crystal graphs that come from
representations of simply laced algebras.

Proposition 2.4.4 of \cite{KMN} states that a crystal with a unique
maximal vertex comes from a representation if and only if it
decomposes as a disjoint union of crystals of representations relative
to the rank 2 subalgebras corresponding to each pair of edge colors.
It therefore suffices to address the problem of locally characterizing
crystal graphs for rank 2 algebras.  The results of \cite{Stem} apply
to the algebras $A_{1}\times A_{1}$ and $A_{2}$; the obvious next case
to consider is $C_{2} \simeq B_{2}$.  In the sequel we call an algebra
doubly laced if all of its regular rank 2 subalgebras are of type
$A_{1}\times A_{1}$, $A_{2}$, or $C_{2}$.  At the end of \cite{Stem},
Stembridge conjectures the following, which we prove in this paper:

\begin{theorem}\label{theorem:int}
    Let $\mathfrak{g}$ be a doubly laced algebra, let $B$ be the
    crystal graph of an irreducible highest weight module of
    $\mathfrak{g}$, and let $v$ be a vertex of $B$ such that $e_{i}v
    \neq 0$ and $e_{j}v \neq 0$, where $e_{i}$ and $e_{j}$ denote two
    different Kashiwara raising operators.  
    Then one of the following is true:
    \begin{enumerate}
        \item $e_{i}e_{j}v = e_{j}e_{i}v$,
    
	\item $e_{i}e_{j}^{2}e_{i}v = e_{j}e_{i}^{2}e_{j}v$ and no
	other sequences of the operators $e_{i}, e_{j}$ with
	length less than or equal to four satisfy such an equality,
    
	\item $e_{i}e_{j}^{3}e_{i}v = e_{j}e_{i}e_{j}e_{i}e_{j}v =
	e_{j}^{2}e_{i}^{2}e_{j}v$ and no other sequences of the
	operators $e_{i}, e_{j}$ with length less than or equal to
	five satisfy such an equality,
    
	\item $e_{i}e_{j}^{3}e_{i}^{2}e_{j}v =
	e_{i}e_{j}^{2}e_{i}e_{j}e_{i}e_{j}v =
	e_{j}e_{i}^{2}e_{j}^{3}e_{i}v = 
	e_{j}e_{i}e_{j}e_{i}e_{j}^{2}e_{i}v$ and 
	no other sequences of the operators $e_{i}, e_{j}$ with length
	less than or equal to seven satisfy such an equality.
    \end{enumerate}
     
    The equivalent statement with $f_{i}$ and $f_{j}$ in place of
    $e_{i}$ and $e_{j}$ also holds.
    
\end{theorem}

We say in these respective cases that $v$ has a degree 2 relation, a
degree 4 relation, a degree 5 relation, or a degree 7 relation above
it.  These may be viewed as combinatorial analogues of the Serre
relations, as observed in \cite{Stem}.  The reader should note that
several of the equalities in the description of degree 5 and degree 7
relations correspond to degree 2 relations within the sequences of
operators.

It suffices to show that Theorem \ref{theorem:int} holds for $C_{2}$
crystals; thanks to the result of \cite{KMN} mentioned above, combined
with the results of \cite{Stem}, the statement automatically extends
to crystals corresponding to representations of any doubly laced
algebra; these algebras are $B_{n}$, $C_{n}$, $F_{4}$, $B_{n}^{(1)}$,
$C_{n}^{(1)}$, $F_{4}^{(1)}$, $A_{n}^{(2)}$, $A_{n}^{(2)\dagger}$,
$D_{n+1}^{(2)}$, and $E_{6}^{(2)}$.

It should be noted that Theorem \ref{theorem:int} does not provide a
local characterization of crystals coming from representations of the
above mentioned algebras.  In order to have such a characterization,
it would be necessary to provide axioms such as those in section 1 of
\cite{Stem} and to show that any graph satisfying those axioms is in
fact a crystal.  Here, we only show the other half of the
characterization; we are assured that any graph with a relation not
explicitly described in the above theorem is in fact not a crystal
over one of these algebras.

\subsection*{Acknowledgments}
The author would like to thank Anne Schilling for suggesting this
problem and for helpful discussions, John Stembridge for very helpful
correspondence on this topic, and Eric Rains and one of the referrees
for the suggestion of working with crystals of type $C_{2}$ rather
than those of type $B_{2}$.

\section{Background on $C_{2}$ crystals}

Recall that a crystal is a colored directed graph in which we
interpret an $i$-colored edge from the vertex $x$ to the vertex $y$ to
mean that $f_{i}x = y$ and $e_{i}y = x$, where $e_{i}$ and $f_{i}$ are
Kashiwara crystal operators.  The vector representation of $C_{2}$ has
the following crystal: \vspace{.2in}
\begin{center}
    \setlength{\unitlength}{1in}
    \begin{picture}(.5,0)
	\put(-1.05,0){\framebox{$1$}}
	\put(-.25,0){\framebox{$2$}}
	\put(.55,0){\framebox{$\bar{2}$}}
	\put(1.3,0){\framebox{$\bar{1}$}}
	\put(-.8,0.05){\vector(1,0){.4}}
	\put(-.6,.16){\footnotesize{1}}
	\put(.05,0.05){\vector(1,0){.4}}
	\put(.27,.16){\footnotesize{2}}
	\put(.8,0.05){\vector(1,0){.4}}
	\put(.95,.16){\footnotesize{1}}
    \end{picture}
    
\end{center}
We realize crystals using tableaux filled with the letters
$\{1,2,\bar{2},\bar{1}\}$ with the total ordering $1 < 2 < \bar{2} <
\bar{1}$.
The following definition is
adapted from that in \cite{KN}.

\begin{definition}
    \label{def:tableau}
	
    A $C_{2}$ Young diagram is a partition with no more than two
    parts; we draw it as a left-justified two-row arrangement of boxes
    such that the second row is no longer than the first.
    
    A $C_{2}$ tableau is a filling of a $C_{2}$ Young diagram by the
    letters of the above alphabet with the following properties:
    \begin{enumerate}
	\item each row is weakly increasing by the ordering in the
	vector representation;

	\item each column is strictly increasing by the ordering in
	the vector representation;

	\item no column may contain $1$ and $\bar{1}$ simultaneously;
	
	\item the configuration 	     
	$\begin{array}{|c|c|}
		\hline
		2 & *  \\
		\hline
		* & \raisebox{-1pt}{$\bar{2}$} \\
		\hline
		 
	\end{array}
	$ does not appear in $T$.

    \end{enumerate}
\end{definition}

\begin{definition}
    \label{def:colword}
    
    Let $T$ be a type $C_{2}$ tableau.  The column word $W$ of $T$ is
    the word on the alphabet $\{ 1, 2, \bar{2}, \bar{1}\}$ consisting
    of $cd$ for each column $
    \begin{array}{|c|} 
	\hline
        d  \\
        \hline
        c  \\
        \hline
         
    \end{array}
    $ in $T$, reading left to right, then followed by each entry
    appearing in a one-row column in $T$, again reading left to right.
    (This could be called the ``reverse far-east reading'', as the
    column word is precisely the reverse of the ``far-east reading''
    used in \cite{HK}).
    
\end{definition}

We now present a definition of the $1$-signature and $2$-signature of
the column word of a type $C_{2}$ tableau, which is easily seen to be
equivalent to the conventional definitions (e.g. \cite{HK}).  Our
definition differs by using the extra symbol $*$ to keep track of
vacant spaces in the signatures.  As in Definition \ref{def:colword}, 
our signatures are in the reverse order from those in \cite{HK}.

\begin{definition}
    \label{def:sigs}
    
    Let $a$ be in the alphabet $\{ 1, 2, \bar{2}, \bar{1}\}$.  Then
    the $1$-signature of $a$ is
    \begin{itemize}
	\item $-$, if $a$ is $\bar{2}$ or $1$;
	
	\item $+$, if $a$ is $\bar{1}$ or $2$;
    
    \end{itemize}
    The $2$-signature of $a$ is
        \begin{itemize}

	\item $-$, if $a$ is $2$;

	\item $+$, if $a$ is $\bar{2}$;

        \item $*$, if $a$ is $1$ or $\bar{1}$.
    \end{itemize}
    
    Let $W$ be the column word of a type $C_{2}$ tableau $T$.  Then
    for $i\in\{1,2\}$ the $i$-signature of $T$ is the word on the
    alphabet $\{ +, -, *\}$ that results from concatenating the
    $i$-signatures of the entries of $W$.

\end{definition}

\begin{definition}
    \label{def:redsig}
    
    Let $S = s_{1}s_{2}\cdots s_{\ell}$ be a signature in the sense of
    Definition \ref{def:sigs}.  The reduced form of $S$ is the word on
    the alphabet $\{+,-,*\}$ that results from iteratively replacing
    every occurance of $+\underbrace{*\cdots*}_{k}-$ in $S$ with
    $\underbrace{*\cdots*}_{k+2}$ until there are no occurances of
    $+\underbrace{*\cdots*}_{k}-$ in $S$.
\end{definition}

The result of applying the Kashiwara operator $e_{i}$ to a tableau $T$
breaks into several cases.  If there are no $+$'s in the reduced form
of the $i$-signature of $T$, we say that $e_{i}T = 0$, where $0$ is a
formal symbol.  Otherwise, let $a$ be the entry corresponding to the
leftmost $+$ in the reduced form of the $i$-signature of $T$.  Then
$e_{i}T$ is the tableau that results from changing $a$ to $e_{i}a$ in
$T$.

Similarly for $f_{i}$, if there are no $-$'s in the reduced form of
the $i$-signature of $T$, we say that $f_{i}T = 0$.  Otherwise, let
$a$ be the entry corresponding to the rightmost $-$ in the reduced
form of the $i$-signature of $T$.  Then $f_{i}T$ is the tableau that
results from changing $a$ to $f_{i}a$ in $T$.

\begin{lemma}
    
    To prove Theorem \ref{theorem:int}, it suffices to prove only the 
    statement regarding $e_{i}$ and $e_{j}$.
 
\end{lemma}

\begin{proof}
    
    For any type $C_{2}$ irreducible highest weight crystal $B$
    corresponding to the module $V$, there is a dual crystal $B^{*}$
    corresponding to the module $V^{*}$.  These crystals are related
    as follows;
    \begin{itemize}
	\item map the highest weight vertex $u_{B}$ of $B$ to the
	lowest weight vertex $\ell_{B^{*}}$ of $B^{*}$,
    
        \item if $v\in B$ is mapped to $w\in B^{*}$, map $f_{i}v$
	to $e_{i}w$.
    \end{itemize}
    
    It is immediate that if the statement regarding $e_{i}$ and
    $e_{j}$ in Theorem \ref{theorem:int} holds, the corresponding
    statement regarding $f_{i}$ and $f_{j}$ holds as well.
    
\end{proof}

\section{Analysis of generic $C_{2}$ tableau}

A generic $C_{2}$ tableau is of the form
\begin{equation}
    \label{eq:gentab}
    \begin{array}{|c|c|c|c|c|c|c|c|c|c|c|c|c|c|c|c|}
	\hline
	\raisebox{-1pt}{$1$} & \cdots & \raisebox{-1pt}{$1$} &
	\raisebox{-1pt}{$1$} & \cdots & \raisebox{-1pt}{$1$}  &
	\raisebox{-1pt}{$2$} & \raisebox{-1pt}{$2$} & \cdots &
	\raisebox{-1pt}{$2$} & 
	\raisebox{-1pt}{$\bar{2}$} & \cdots &
	\raisebox{-1pt}{$\bar{2}$} & \raisebox{-1pt}{$\bar{1}$} & \cdots &
	\raisebox{-1pt}{$\bar{1}$}  \\
	\cline{1-16}
	\raisebox{-1pt}{$2$} & \cdots & \raisebox{-1pt}{$2$} &
	\raisebox{-1pt}{$\bar{2}$} & \cdots &
	\raisebox{-1pt}{$\bar{2}$} &
	\raisebox{-1pt}{$\bar{2}$} & \raisebox{-1pt}{$\bar{1}$} & \cdots &
	\raisebox{-1pt}{$\bar{1}$} & \raisebox{-1pt}{$\bar{1}$} & \cdots &
	\raisebox{-1pt}{$\bar{1}$} 
	& \multicolumn{3}{}{c}  \\
	\cline{1-13}
    \end{array}
\end{equation}
where 
\begin{itemize}
    \item any column may be omitted;

    \item  any of the columns other than
    $\begin{array}{|c|} \hline
	2 \\
	\hline
	\raisebox{-1pt}{$\bar{2}$}  \\
	\hline
	 
    \end{array}
    $ may be repeated an arbitrary number of times;

    \item the bottom row
    may be truncated at any point.

\end{itemize}
 
We are interested in how the Kashiwara operators $e_{1}$ and $e_{2}$
act on this tableaux, so we must determine where the left-most $+$
appears in the reduced form of the signatures of the tableau.  The
relevant $+$'s in the signatures of a generic tableau naturally fall
into two groups as described by definition \ref{def:blocks}.

\begin{definition}
    \label{def:blocks}
    Let $T$ be a $C_{2}$ tableau.  
    \begin{itemize}
	\item We define the left block of $+$'s in the $1$-signature
	of $T$ to be those $+$'s from $2$'s in the top row and
	$\bar{1}$'s in the bottom row.  If no such entries appear in
	$T$, we say that the left block of $+$'s in the $1$-signature
	of $T$ has size $0$ and its left edge is located on the
	immediate left of symbol coming from the leftmost $\bar{2}$ or
	$\bar{1}$ in the top row of $T$.  If there is furthermore no
	such entry, its left edge is located at the right end of the
	$1$-signature of $T$.
    
	\item We define the right block of $+$'s in the $1$-signature
	of $T$ to be those $+$'s from $\bar{1}$ in the top row of $T$.
	If no such entry appears in $T$, we say that the right block
	of $+$'s in the $1$-signature of $T$ has size $0$ and its left
	edge is located at the right end of the $1$-signature of $T$.
    
	\item We define the left block of $+$'s in the $2$-signature
	of $T$ to be those $+$'s from $\bar{2}$ in the bottom row.  If
	such an entry does not appear in $T$, we say that the left
	block of $+$'s in the $2$-signature of $T$ has size $0$ and
	its left edge is located on the immediate left of the $*$ in
	the $2$-signature coming from the leftmost $\bar{1}$ in the
	bottom row of $T$.  If $T$ has no $\bar{1}$'s in the bottom
	row, we say that its left edge is located at the right end of
	the $2$-signature of $T$.
    
	\item We define the right block of $+$'s in the $2$-signature
	of $T$ to be those $+$'s from $\bar{2}$ in the top row of $T$.
	If such an entry does not appear in $T$, we say that the right
	block of $+$'s in the $2$-signature of $T$ has size $0$ and
	its left edge is located on the immediate left of the $*$ in
	the $2$-signature coming from the leftmost $\bar{1}$ in the
	top row of $T$.  If there are furthermore no $\bar{1}$'s in
	the top row of $T$, we say that its left edge is located at
	the right end of the $2$-signature of $T$.
	
    \end{itemize}
    
    In the above cases when a block of $+$'s has positive size we say
    that its left edge is on the immediate left of its leftmost $+$.
    
\end{definition}

Motivated by this definition, we define the following statistics on a
$C_{2}$ tableaux $T$.
\begin{itemize}
    \item $A(T)$ is the number of $\bar{2}$'S in the top row of $T$,

    \item $B(T)$ is the number of $2$'s in the top row of $T$ plus the
    number of $\bar{1}$'s in the bottom row of $T$,

    \item $C(T)$ is the number of $2$'s in the top row of $T$,

    \item $D(T)$ is the number of $\bar{2}$'s in the bottom row of
    $T$.
\end{itemize}

\begin{example}
    Let 
    $$
    T = 
    \begin{array}{|c|c|c|c|c|c|c}
        
	\cline{1-6}
        \raisebox{-1pt}{1} & \raisebox{-1pt}{1} & \raisebox{-1pt}{2} & \raisebox{-1pt}{$\bar{2}$} & \raisebox{-1pt}{$\bar{2}$} & \raisebox{-1pt}{$\bar{1}$} &   \\
        \cline{1-6}
        \raisebox{-1pt}{2} & \raisebox{-1pt}{$\bar{2}$} & \raisebox{-1pt}{$\bar{2}$} & \raisebox{-1pt}{$\bar{1}$} & \multicolumn{3}{c}{}    \\
        
        \cline{1-4} 
    \end{array}
    .
    $$
    Then $A(T) = 2$, $B(T) = 2$, $C(T)  = 1$, and $D(T) = 2$.
\end{example}

\begin{claim} \label{claim:lr}
    $\begin{array}{c}
        
          \\

    \end{array}
    $
    \begin{itemize}
	\item If $e_{1}$ acts on a tableau $T$ with $A(T) < B(T)$, the
	entry on which $e_{1}$ acts corresponds to a symbol in the
	left block of $+$'s in the $1$-signature of $T$;
    
	\item If $e_{1}$ acts on a tableau $T$ with $A(T) \geq B(T)$,
	the entry on which $e_{1}$ acts corresponds to a symbol in
	the right block of $+$'s in the $1$-signature of $T$;
    
	\item If $e_{2}$ acts on a tableau $T$ with $C(T) < D(T)$, the
	entry on which $e_{2}$ acts corresponds to a symbol in the
	left block of $+$'s in the $2$-signature of $T$;
    
	\item If $e_{2}$ acts on a tableau $T$ with $C(T) \geq D(T)$,
	the entry on which $e_{2}$ acts corresponds to a symbol in
	the right block of $+$'s in the $2$-signature of $T$.
    \end{itemize}
    
\end{claim}

We now show that the entries in $T$ on which a sequence
$e_{i_{1}}\cdots e_{i_{\ell}}$ acts are determined by which blocks of
$+$'s correspond to those entries.  To achieve this, we verify that
the left edge of a block of $+$'s can be changed only by acting on
that block of $+$'s.  This goal motivates the following notation.

We write $e_{1}^{{\bf \ell}}$ to indicate the Kashiwara operator
$e_{1}$ when applied to a tableau $T$ such that $A(T) < B(T)$ and
$e_{1}^{{\bf r}}$ to indicate the Kashiwara operator $e_{1}$ when
applied to a tableau $T$ such that $A(T) \geq B(T)$.  Similarly, we
write $e_{2}^{{\bf \ell}}$ to indicate the Kashiwara operator $e_{2}$
when applied to a tableau $T$ such that $C(T) < D(T)$ and $e_{2}^{{\bf
r}}$ to indicate the Kashiwara operator $e_{2}$ when applied to a
tableau $T$ such that $C(T) \geq D(T)$.  Note that these are not new
operators: we simply use the superscript notation to record additional
information about how the operators act on specific tableaux.

\begin{claim} \label{prop:leftedge}
    Let $T$ be a tableau such that $e_{1}^{{\bf r}}T \neq 0$.  Then
    the left edges of both the left and right blocks of $+$'s in the
    $2$-signature of $e_{1}T$ and the left edge of the left block of
    $+$'s in the $1$-signature of $e_{1}T$ are in the same place as
    they are in the signatures of $T$, and the left edge of the right
    block of $+$'s in the $1$-signature of $e_{1}T$ is one position to
    the right of that in $T$.  
    
    Symmetrically, if $e_{1}^{\ell}T \neq 0$, the left edges of the
    blocks in the $2$-signature and the right block in the
    $1$-signature are unchanged and the left edge of the left block in
    the $1$-signature moves one position to the right; if $e_{2}^{{\bf
    r}}T \neq 0$, the left edges of the blocks in the $1$-signature
    and the left block in the $2$-signature are unchanged and the left
    edge of the right block in the $2$-signature moves one position to
    the right; and if $e_{2}^{\ell}T \neq 0$, the left edges of the
    blocks in the $1$-signature and the right block in the
    $2$-signature are unchanged and the left edge of the left block in
    the $2$-signature moves one position to the right.
\end{claim}

\begin{proof}
    Let $i = 1$ and $j = 2$ or vice versa, and let ${\bf x} = {\bf
    \ell}$ and ${\bf y} = {\bf r}$ or vice versa.  It is clear from
    the combinatorially defined action of $e_{i}$ on $C_{2}$ tableaux
    that if $e_{i}^{{\bf x}}T \neq 0$, the left edge of the {\bf y}
    block of $+$'s in the $i$-signature of $e_{i}T$ is in the same
    place as in the $i$-signature of $T$, and that the left edge of
    the {\bf x} block of $+$'s in the $i$-signature of $e_{i}T$ is one
    space to the right of its position in the $i$-signature of $T$.
    We may therefore devote our attention to the $j$-signature in each
    of the four cases of concern.

    First, consider the case of $e_{1}^{{\bf \ell}}T\neq 0$.  By Claim
    \ref{claim:lr}, we know that this operator changes a $2$ in the
    top row to a $1$ or a $\bar{1}$ in the bottom row to a $\bar{2}$.
    In the first case, no $+$'s are added to the $2$-signature, and no
    change is made to those entries of interest to the location of a
    block of $+$'s of size $0$ in the $2$-signature.  In the other
    case, one $+$ is added to the left block of $+$'s in the
    $2$-signature; this addition is to the right of the left edge of
    this block.  Finally, the right block of $+$'s in the
    $2$-signature of $e_{1}T$ is the same as in the $2$-signature of
    $T$ in any case.

    Next, consider the case of $e_{1}^{{\bf r}}T\neq 0$.  By Claim
    \ref{claim:lr}, we know that this operator changes a $\bar{1}$ in the
    top row of $T$ into a $\bar{2}$.  This adds one $+$ to the right block
    of $+$'s in the $2$-signature to the right of its left edge and makes
    no change to the left block of $+$'s.

    Now, consider the case of $e_{2}^{{\bf \ell}}T\neq 0$.  By Claim
    \ref{claim:lr}, we know that this operator changes a $\bar{2}$ in the
    bottom row to a $2$.  This does not contribute a $+$ to either the
    left or right blocks of the $1$-signature of $e_{2}T$, nor does it
    pertain to the location of a block of $+$'s of size $0$ in the
    $1$-signature, so the left edges in this signature are the same as in
    the $1$-signature of $T$.

    Finally, we consider the case of $e_{2}^{{\bf r}}T\neq 0$.  By Claim
    \ref{claim:lr}, we know that this operator changes a $\bar{2}$ in the
    top row to a $2$.  This has the effect of adding a $+$ to the left
    block of $+$'s in the $1$-signature to the right of the left edge of
    this block.  Finally, the right block of $+$'s in the $1$-signature of
    $e_{2}T$ is the same as in the $1$-signature of $T$.

\end{proof}

\begin{corollary}
    Let $T$ be a tableau such that $e_{1}^{\ell}T \neq 0$, and let $E$
    be a sequence of operators from the set $\{ e_{1}^{{\bf r}},
    e_{2}^{\ell}, e_{2}^{{\bf r}} \}$ such that $ET \neq 0$.  Then
    $e_{1}^{\ell}$ acts on the same entry in $T$ as it does in $ET$.
    The symmetric statements corresponding to the cases of Claim
    \ref{prop:leftedge} hold as well.
\end{corollary}

The following four Sublemmas state that the relative values of $A(T),
B(T), C(T),$ and $D(T)$ not only determine where $e_{i}$ acts within a
tableau, but also what the values of $A(e_{i}T), B(e_{i}T),
C(e_{i}T),$ and $D(e_{i}T)$ are.  This will be an invaluable tool for
our analysis in section \ref{sec:main_proof}.

\begin{sublemma} \label{sublemma:alb}
    Suppose $T$ is a tableau such that $e_{1}$ acts on the left block
    of $+$'s in the $1$-signature of $T$ (i.e., such that $A(T) <
    B(T)$).  Then $A(e_{1}T) = A(T)$, $B(e_{1}T) = B(T) - 1$, and
    $C(e_{1}T) - D(e_{1}T) = C(T) - D(T) - 1$.
\end{sublemma}

\begin{proof}
    
We have two cases to consider; $e_{1}$ may act by changing a $2$ to a
$1$ in the top row or a $\bar{1}$ to a $\bar{2}$ in the bottom row.
In both of these cases, it is easy to see that the number of
$\bar{2}$'s in the top row is unchanged and the number of $2$'s in the
top row plus the number of $\bar{1}$'s in the bottom row is diminished
by one; hence $A(e_{1}T) = A(T)$ and $B(e_{1}T) = B(T) - 1$.

Observe that in the case of a $\bar{1}$ changing into a $\bar{2}$ in
the bottom row, the content of the top row is unchanged, but the
number of $\bar{2}$'s in the bottom row is increased by $1$.  In the
case of a $2$ changing into a $1$ in the top row, the bottom row is
unchanged, but the number of $2$'s in the top row is decreased by $1$.
In both of these cases, we find that $C(e_{1}T) - D(e_{1}T) = C(T) -
D(T) - 1$.

\end{proof}

\begin{sublemma} \label{sublemma:agb}
    Suppose $T$ is a tableau such that $e_{1}$ acts on the right block
    of $+$'s in the $1$-signature of $T$ (i.e., such that $A(T) \geq
    B(T)$).  Then $A(e_{1}T) = A(T) + 1$, $B(e_{1}T) = B(T)$,
    $C(e_{1}T) = C(T)$, and $D(e_{1}T) = D(T)$.
\end{sublemma}

\begin{proof}

Since the right block of $+$'s in the $1$-signature comes entirely
from $\bar{1}$'s in the top row of $T$, it follows that acting by
$e_{1}$ changes one of these $\bar{1}$'s into a $\bar{2}$.  We
immediately see that the number of $\bar{2}$'s in the top row
increases by $1$, and that the number of $2$'s in the top row and
$\bar{2}$'s and $\bar{1}$'s in the bottom row are all unchanged.

\end{proof}

\begin{sublemma} \label{sublemma:cld}
    Suppose $T$ is a tableau such that $e_{2}$ acts on the left block
    of $+$'s in the $1$-signature of $T$ (i.e., such that $C(T) <
    D(T)$).  Then $A(e_{2}T) = A(T)$, $B(e_{2}T) = B(T)$, $C(e_{2}T) =
    C(T)$, and $D(e_{2}T) = D(T) - 1$.
\end{sublemma}

\begin{proof}

The entry on which $e_{2}$ acts is a $\bar{2}$ the bottom row, which
will be changed into a $2$.  We immediately see that the number of
$\bar{2}$'s in the bottom row decreases by $1$, and that the number of
$2$'s and $\bar{2}$'s in the top row and $\bar{1}$'s in the bottom row
are all unchanged.

\end{proof}

\begin{sublemma} \label{sublemma:cgd}
    Suppose $T$ is a tableau such that $e_{2}$ acts on the right block
    of $+$'s in the $1$-signature of $T$ (i.e., such that $C(T) \geq
    D(T)$).  Then $A(e_{2}T) = A(T) - 1$, $B(e_{2}T) = B(T) + 1$,
    $C(e_{2}T) = C(T) + 1$, and $D(e_{2}T) = D(T)$.
\end{sublemma}

\begin{proof}

In this case $e_{2}$ will change a $\bar{2}$ to a $2$ in the top row.
It is easy to see that the number of $2$'s in the top row increases by
$1$ and the number of $\bar{2}$'s in the bottom row is unchanged;
hence $C(e_{2}T) = C(T) + 1$ and $D(e_{2}T) = D(T)$.
  
Likewise, since the number of $\bar{2}$'s in the top row is decreased
by $1$ and the number of $2$'s in the top row is increased by $1$, we
find that $A(e_{2}T) = A(T) - 1$ and $B(e_{2}T) = B(T) + 1$.

\end{proof}

\section{Proof of Theorem \ref{theorem:int}}
\label{sec:main_proof}

We are now equipped to begin addressing Theorem
\ref{theorem:int}.  It is proved as a consequence of Lemmas
\ref{lemma:deg2:1} through \ref{lemma:deg7}, each of which deals with
a certain case of the relative values of $A(T), B(T), C(T),$ and
$D(T)$.  To see that these cases are exhaustive, refer to Table
\ref{relstats}.

\begin{table}
    $
    \begin{array}{c|c|c|c|c|}
	 & A<B & A=B & A=B+1 & A>B+1  \\
	\hline
	 C<D & 2 & 2 & 2 & 2  \\
	\hline
	 C=D & 4 & 7 & 4 & 2  \\
	\hline
	 C>D & 2 & 5 & 4 & 2  \\
	\hline
	 
    \end{array}
    $

    \caption{Degree of relation over $T$, given its $ABCD$ statistics} 
    \protect\label{relstats}
\end{table}


\begin{lemma}
    \label{lemma:deg2:1}
    Suppose $T$ is a tableau such that $C(T) < D(T)$, $e_{1}T\neq 0$, and
    $e_{2}T\neq 0$.  Then $T$ has a degree 2 relation above it.
\end{lemma}

\begin{proof}

    From Claim \ref{claim:lr}, we know that $e_{2}$ acts on the left
    block of $+$'s, and by Sublemma \ref{sublemma:cld}, we know that
    $A(e_{2}T) = A(T)$ and $B(e_{2}T) = B(T)$; it follows that
    $e_{1}e_{2}T \neq 0$, and that $e_{1}$ acts on the same entry in
    $e_{2}T$ as it does in $T$.  Furthermore, by Sublemmas
    \ref{sublemma:alb} and \ref{sublemma:agb}, we know that either
    $C(e_{1}T) - D(e_{1}T) = C(T) - D(T) - 1$ or $C(e_{1}T) = C(T)$
    and $D(e_{1}T) = D(T)$; in either case, we still find that
    $C(e_{1}T) < D(e_{1}T)$.  Since $C(e_{1}T) \geq 0$, we are assured
    that $D(e_{1}T) \geq 1$, and thus $e_{2}e_{1}T \neq 0$.  We
    conclude that $e_{2}$ acts on the same entry in $e_{1}T$ as it
    does in $T$.
    
\end{proof}


\begin{lemma}
    \label{lemma:deg2:2}
    Suppose $T$ is a tableau such that $A(T) > B(T) + 1$, $e_{1}T\neq
    0$, and $e_{2}T\neq 0$.  Then $T$ has a degree 2 relation above
    it.
\end{lemma}

\begin{proof}
    
    From Claim \ref{claim:lr}, we know that $e_{1}$ acts on the right
    block of $+$'s, and by Sublemma \ref{sublemma:agb}, we know that
    $C(e_{1}T) = C(T)$ and $D(e_{1}T) = D(T)$; thus $e_{2}e_{1}T \neq
    0$, and $e_{2}$ acts on the same entry in $e_{1}T$ as it does in
    $T$.  Furthermore, by Sublemmas \ref{sublemma:cld} and
    \ref{sublemma:cgd}, we know that either $A(e_{2}T) = A(T)$ and
    $B(e_{2}T) = B(T)$ or $A(e_{2}T) = A(T) - 1$ and $B(e_{2}T) = B(T)
    + 1$; in either case, we find that $A(e_{2}T) > B(e_{2}T)$ and the
    size of the right block of $+$'s in the $1$-signature is not
    diminished.  We therefore conclude that $e_{1}e_{2}T \neq 0$ and
    that $e_{1}$ acts on the same entry in $e_{2}T$ as it does in
    $T$.

\end{proof}


\begin{lemma}\label{lemma:deg2:3}
    Suppose $T$ is a tableau such that $A(T) < B(T)$, $C(T) > D(T)$,
    $e_{1}T\neq 0$, and $e_{2}T\neq 0$.  Then $T$ has a degree 2
    relation above it.
\end{lemma}

\begin{proof}
    
    By Claim \ref{claim:lr}, we know that $e_{1}$ acts on the left
    block of $+$'s in $T$ and $e_{2}$ acts on the right block of $+$'s
    in $T$.  By Sublemma \ref{sublemma:cgd}, we know that $A(e_{2}T) =
    A(T) - 1$ and $B(e_{2}T) = B(T) + 1$.  It follows that $A(e_{2}T)
    < B(e_{2}T)$, and since $A(e_{2}T) \geq 0$, this ensures that
    $B(e_{2}T) \geq 1$, and thus $e_{1}e_{2}T \neq 0$.  We conclude
    that $e_{1}$ acts on the same entry in $e_{2}T$ as it does in
    $T$.  Furthermore, by Sublemma \ref{sublemma:alb}, we know that
    $C(e_{1}T) - D(e_{1}T) = C(T) - D(T) - 1$; it follows that
    $C(e_{1}T) \geq D(e_{1}T)$.  Since we also know that the size of
    the right block of $+$'s in the $2$-signature of $e_{1}T$ is at
    least as large as that of $T$, it is the case that $e_{2}e_{1}T
    \neq 0$, and so $e_{2}$ acts on the same entry in $e_{1}T$ as it
    does in $T$.

\end{proof}


To prove Lemmas \ref{lemma:deg4:1} through \ref{lemma:deg7}, we must
not only show that the given sequences of operators act on the same
entries, but also that no pair of homogeneous sequences of operators
(i.e., a pair $(P_{1},P_{2})$ such that $P_{1}$ and $P_{2}$ have the
same number of instances of $e_{1}$ and $e_{2}$) with shorter or equal
length act on the same entries.  To assist in our illustration of this
fact, we will refer to figures that encode the generic behavior of all
sequences of operators on a tableau with content as specified by the
hypothesis of each lemma.  Table \ref{piclegend} is a legend for the
figures used to prove Lemmas \ref{lemma:deg4:1} through
\ref{lemma:deg5}.  In the picture used to prove Lemma
\ref{lemma:deg7}, we instead use an edge pointing down to indicate
acting by $e_{1}$ and an edge pointing up to indicate acting by
$e_{2}$; otherwise the legend is the same.

To assist in proving that the sequences in question do not kill our
tableaux, we have the following Sublemma.

\begin{sublemma}
    \label{sublemma:dashed}
    Let $E$ be a dashed edge from $v$ up to $w$; i.e., an operator
    $e_{i}^{\ell}$ acts on $v$ to produce $w$.  Then $e_{i}v \neq 0$.
\end{sublemma}

\begin{proof}
    
    The Kashiwara operator $e_{i}$ acts on the left block of $+$'s of
    a tableau $T$ precisely when $A(T) < B(T)$ or $C(T) < D(T)$ in the
    cases of $i = 1$ or $i = 2$, respectively.  Since these numbers
    are all non-negative integers, we conclude that $B(T) > 0$ or
    $D(T) > 0$.  Since these statistics indicate the number of $+$'s
    in the left block of their respective signatures, we are assured
    that there is an entry on which $e_{i}$ can act.
    
\end{proof}

Thus it suffices to prove that the solid edges in the paths of concern
do not produce $0$.

\begin{table}
    
    $
    \begin{array}{|c|c|}
	\hline
	\textrm{edge pointing to the right} & \textrm{acting by }e_{1}  \\
	\hline
	\textrm{edge pointing to the left} & \textrm{acting by }e_{2}  \\
	\hline
	\textrm{solid edge} & \textrm{acting on the right block of }+\textrm{'s}  \\
	\hline
	\textrm{dashed edge} & \textrm{acting on the left block of }+\textrm{'s}  \\
	\hline
	\textrm{vertex labeled by }(d_{1},d_{2}) & \textrm{tableau
	}T\textrm{ with statistics such that }\\&A(T) = B(T) +
	d_{1}\textrm{ and }C(T) = D(T) + d_{2} \\
	\hline
	 
    \end{array}
    $

    \caption{Legend for Figures \ref{piclem4} through \ref{piclem7}}
    \protect\label{piclegend}
\end{table}


\begin{lemma}
    \label{lemma:deg4:1}
    Suppose $T$ is a tableau such that $A(T) = B(T) + 1$, $C(T) \geq
    D(T)$, $e_{1}T\neq 0$, and $e_{2}T\neq 0$.  Then $T$ has a degree
    4 relation above it.
\end{lemma}

\begin{figure}
    \centerline{\includegraphics[height=2.5in]{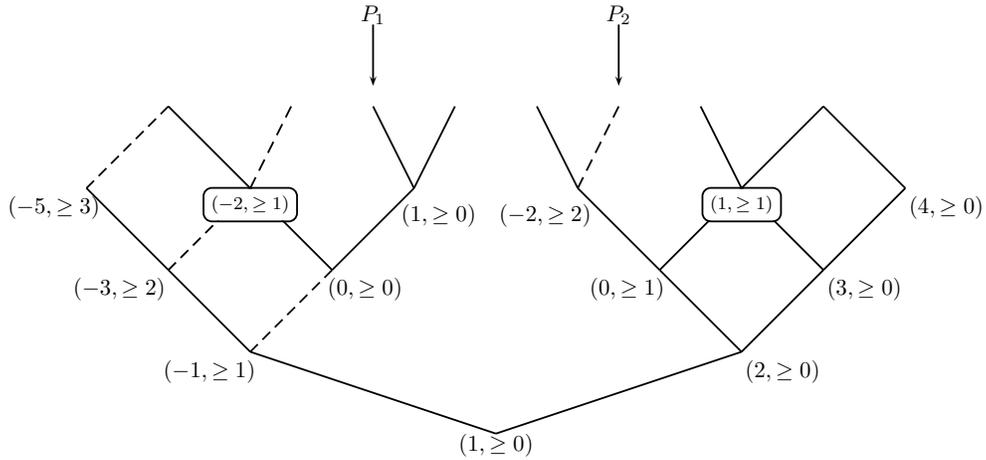}}
    \caption{Picture for Lemma \ref{lemma:deg4:1}}
    \protect\label{piclem4}
\end{figure}

\begin{proof}

    We must first confirm that the sequences $e_{1}e_{2}^{2}e_{1}$ and
    $e_{2}e_{1}^{2}e_{2}$ do not produce $0$ when applied to $T$.
    First, observe that since $e_{1}$ acts on the right block of $+$'s
    in $T$, it changes a $\bar{1}$ to a $\bar{2}$ in the top row.
    This adds a $+$ to the reduced $2$-signature of the tableaux, so
    we know that $e_{2}^{2}e_{1}T \neq 0$.  By Sublemma
    \ref{sublemma:dashed}, we know that $e_{1}e_{2}^{2}e_{1}T \neq 0$.
    On the other hand, we know that $e_{2}$ acts on $T$ by changing a
    $\bar{2}$ to a $2$ in the top row; this means that the reduced
    $1$-signature of $e_{2}T$ has a single $+$ in the left block and
    its right block has at least one $+$, as did the $1$-signature of
    $T$.  We conclude that $e_{1}^{2}e_{2}T \neq 0$.  We know that in
    $e_{1}e_{2}T$, $e_{1}$ changes a $\bar{1}$ to a $\bar{2}$ in the
    top row; since the $+$'s in the $2$-signature from this entry
    cannot be paired with any $-$'s, we conclude that
    $e_{2}e_{1}^{2}e_{2}T \neq 0$.

    Now that we know that neither of these sequences produces $0$ when
    applied to $T$, it is clear that we have $e_{1}e_{2}^{2}e_{1}T =
    e_{2}e_{1}^{2}e_{2}T$, as the paths $P_{1}$ and $P_{2}$ leading
    from the base of the graph in Figure \ref{piclem4} to the
    indicated leaves both have one solid right edge, one dashed right
    edge, and two solid left edges.  We must now confirm that among
    all pairs $(Q_{1},Q_{2})$ of increasing paths from the base in
    these graphs such that $Q_{1}$ begins by following the left edge
    and $Q_{2}$ begins by following the right edge, $(P_{1}, P_{2})$
    is the only pair with the same number of each type of edge.

    Since the right edge from the base of the graph is solid, our
    candidate for $Q_{1}$ must have a solid right edge.  Inspecting
    the graph tells us that this path must begin with the path
    corresponding to $e_{1}^{2}e_{2}$.  This path has a dashed left
    edge, and the only candidate for $Q_{2}$ with this feature is in
    fact $P_{2}$, which has two solid right edges.  The only way to
    extend $e_{1}^{2}e_{2}$ to have the same edge content as $P_{2}$
    is by extending it to $P_{1}$.

\end{proof}


\begin{lemma}
    \label{lemma:deg4:2}
    Suppose $T$ is a tableau such that $A(T) < B(T)$, $C(T) = D(T)$,
    $e_{1}T\neq 0$, and $e_{2}T\neq 0$.  Then $T$ has a degree 4
    relation above it.
\end{lemma}

\begin{figure}
    \centerline{\includegraphics[height=2.5in]{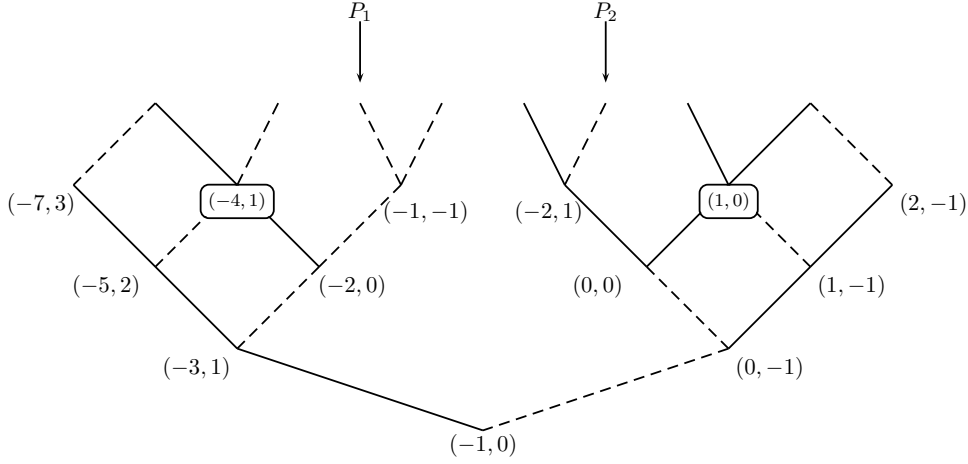}}
    \caption{Picture 1 for Lemma \ref{lemma:deg4:2}}
    \protect\label{pic1lem5}
\end{figure}

\begin{figure}
    \centerline{\includegraphics[height=2.5in]{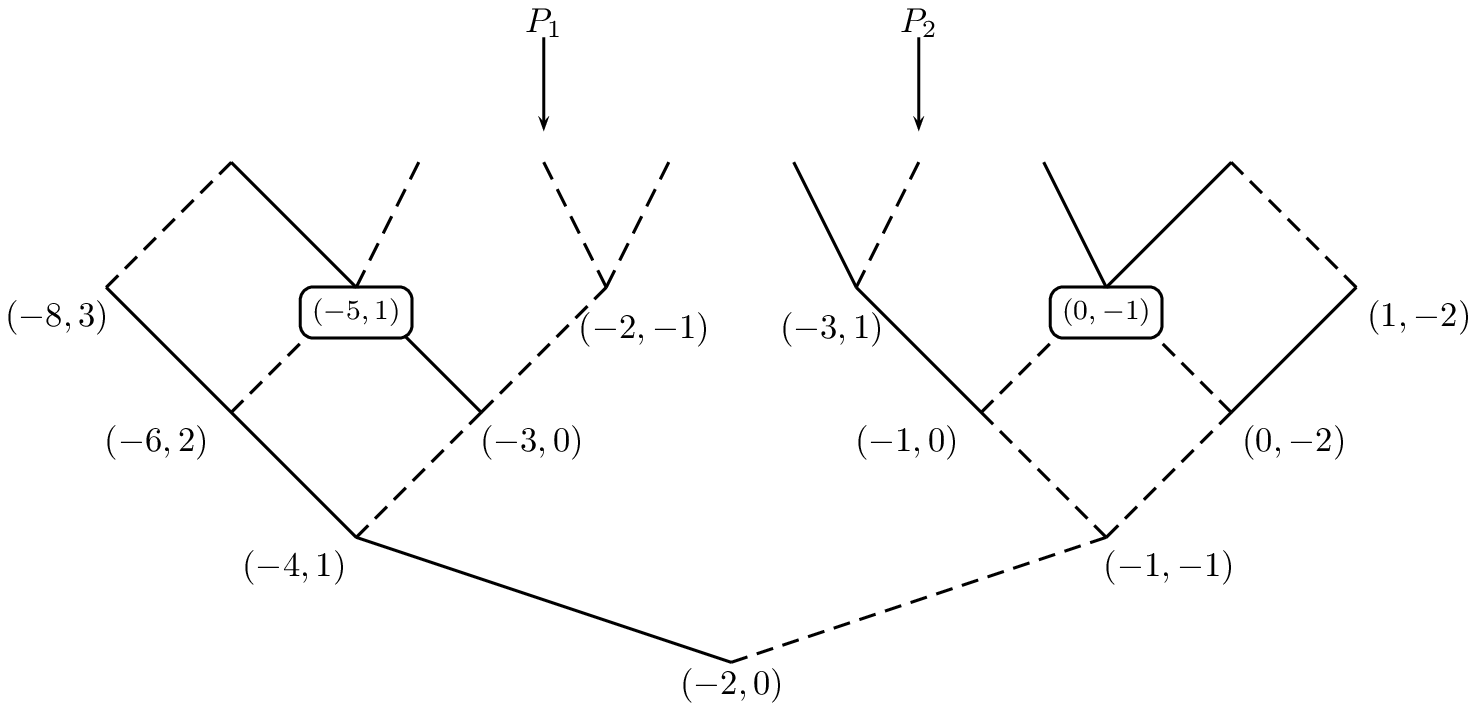}}
    \caption{Picture 2 for Lemma \ref{lemma:deg4:2}}
    \protect\label{pic2lem5}
\end{figure}

\begin{figure}
    \centerline{\includegraphics[height=2.5in]{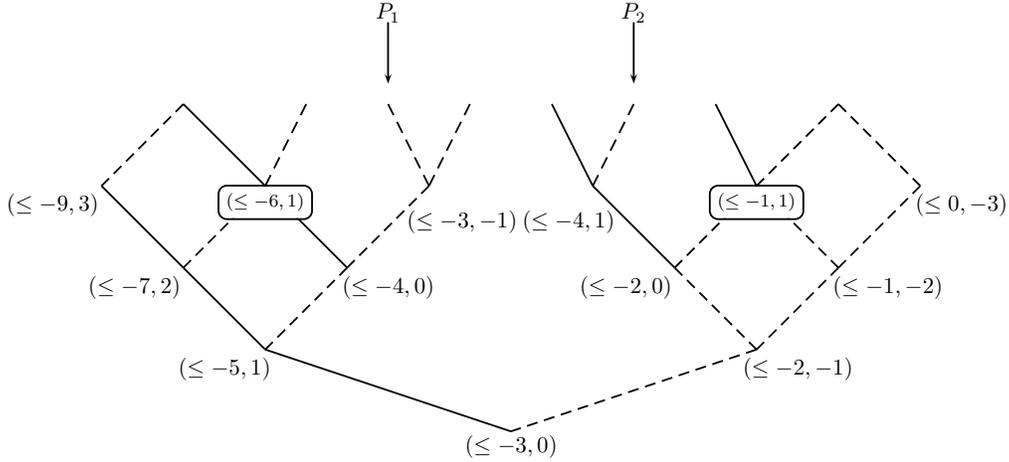}}
    \caption{Picture 3 for Lemma \ref{lemma:deg4:2}}
    \protect\label{pic3lem5}
\end{figure}

\begin{proof}

    We must first confirm that the sequences $e_{1}e_{2}^{2}e_{1}$ and
    $e_{2}e_{1}^{2}e_{2}$ do not produce $0$ when applied to $T$.  By
    Sublemma \ref{sublemma:dashed} it suffices to show that
    $e_{2}^{2}e_{1}T \neq 0$, since $e_{2}T \neq 0$ by assumption and
    all other edges in $P_{1}$ and $P_{2}$ are dashed.  To see this,
    simply observe that there is at least one $+$ in the right block
    of the reduced $2$-signature of $T$; the sequence $e_{2}e_{1}$
    acts on the left blocks of $+$'s, so the corresponding entry
    remains available for $e_{2}$ to act on.

    Now that we know that neither of these sequences produces $0$ when
    applied to $T$, it is clear that we have $e_{1}e_{2}^{2}e_{1}T =
    e_{2}e_{1}^{2}e_{2}T$, as the paths $P_{1}$ and $P_{2}$ leading
    from the base of the graph in Figures \ref{pic1lem5} through
    \ref{pic3lem5} to the indicated leaves have two dashed right
    edges, one solid left edge, and one dashed left edge.  We must now
    confirm that among all pairs $(Q_{1},Q_{2})$ of increasing paths
    from the base in these graphs such that $Q_{1}$ begins by
    following the left edge and $Q_{2}$ begins by following the right
    edge, $(P_{1}, P_{2})$ is the only pair with the same number of
    each type of edge.

    This is easy to see by the following argument.  Every candidate
    for $Q_{2}$ (i.e., every path in the right half of the graphs in
    Figures \ref{pic1lem5} through \ref{pic3lem5}) has at least one
    dashed left edge.  The only candidate for $Q_{1}$ (i.e., the only
    path in the left half of the graphs in Figures \ref{pic1lem5}
    through \ref{pic3lem5}) with a dashed left edge is $P_{1}$.  By
    inspecting Figures \ref{pic1lem5} through \ref{pic3lem5}, $P_{2}$
    is the only candidate for $Q_{2}$ with two dashed right edges and
    one solid left edge.

\end{proof}


\begin{lemma}
    \label{lemma:deg5}
    Suppose $T$ is a tableau such that $A(T) = B(T)$, $C(T) \geq D(T)
    + 1$, $e_{1}T\neq 0$, and $e_{2}T\neq 0$.  Then $T$ has a degree 5
    relation above it.
\end{lemma}

\begin{figure}
    \centerline{\includegraphics[height=2.5in]{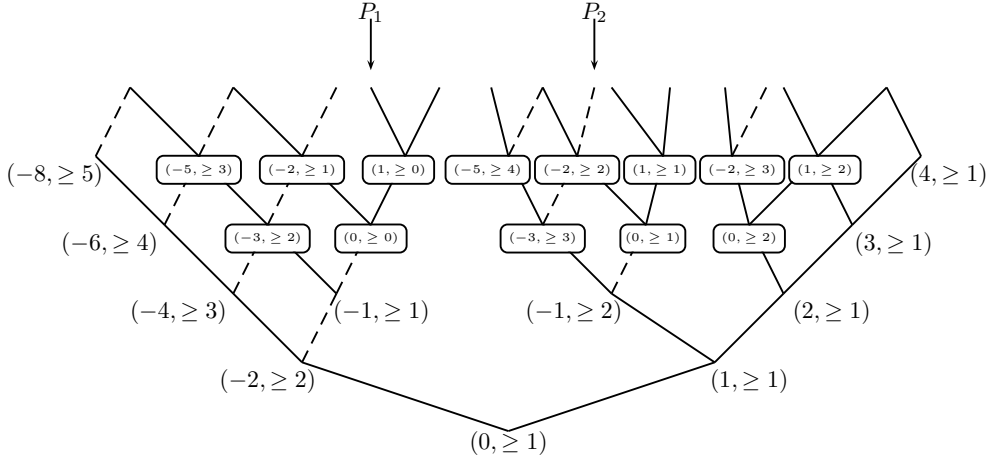}}
    \caption{Picture for Lemma \ref{lemma:deg5}}
    \protect\label{pic1lem6}
\end{figure}

\begin{proof}

    We must first confirm that the sequences $e_{2}e_{1}^{3}e_{2}$ and
    $e_{1}e_{2}e_{1}e_{2}e_{1}$ do not produce $0$ when applied to
    $T$.  First, note that there is at least one $+$ in the right
    block of the reduced $1$-signature of $T$.  Since $e_{2}$ acts on
    the right block of $+$'s in the $2$-signature of $T$, there are as
    many $+$'s in the right block of the $1$-signature of $e_{2}T$ as
    in that of $e_{1}T$.  Observe that $A(e_{2}T) = B(e_{2}T) - 2$, so
    we know that there are additionally two $+$'s in the left block of
    the reduced $1$-signature of $e_{2}T$.  This implies that
    $e_{1}^{3}e_{2}T \neq 0$.  The third of these applications of
    $e_{1}$ changes a $\bar{1}$ to a $\bar{2}$ in the top row; the $+$
    in the $2$-signature of $e_{1}^{3}e_{2}T$ entry cannot be
    bracketed, so we know that $e_{2}e_{1}^{3}e_{2}T \neq 0$.  On the
    other hand, we know that the right block of the reduced
    $2$-signature of $T$ has at least one $+$.  Since $e_{1}$ changes
    a $\bar{1}$ to a $\bar{2}$ in the top row of $T$, we know that
    $e_{2}$ will change this entry to a $2$ so that the right block
    of the reduced signature of $e_{1}T$ has at least two $+$'s that
    cannot be bracketed by $-$'s.  The leftmost of these entries will
    be acted upon by $e_{2}$, so $e_{2}e_{1}T \neq 0$.  Furthermore,
    since $A(e_{2}e_{1}T) = B(e_{2}e_{1}T) - 1$, we know that
    $e_{1}e_{2}e_{1}T \neq 0$.  At least one $+$ remains in the right
    block of the reduced $2$-signature of $e_{1}e_{2}e_{1}T$, so
    $e_{2}e_{1}e_{2}e_{1}T \neq 0$.  Finally, since
    $A(e_{2}e_{1}e_{2}e_{1}T) = B(e_{2}e_{1}e_{2}e_{1}T) - 2$, we know
    that $e_{1}e_{2}e_{1}e_{2}e_{1}T \neq 0$.
    
    Now that we know that neither of these sequences produces $0$ when
    applied to $T$, it is clear that we have $e_{2}e_{1}^{3}e_{2}T =
    e_{1}e_{2}e_{1}e_{2}e_{1}T$, since the paths $P_{1}$ and $P_{2}$
    leading from the base of the graph in Figure \ref{pic1lem6} to the
    indicated leaves have no solid left edges, two dashed left edges,
    one solid right edge, and two dashed right edges.  Note that these
    paths are equivalent to $e_{1}^{2}e_{2}^{2}e_{1}T$, due to the
    degree 2 relation above $e_{2}e_{1}T$; we may denote this
    alternative path by $P'_{2}$.  We must now confirm that among all
    pairs $(Q_{1},Q_{2})$ of increasing paths from the base in these
    graphs such that $Q_{1}$ begins by following the left edge and
    $Q_{2}$ begins by following the right edge, $(P_{1}, P_{2})$ and
    $(P_{1},P'_{2})$ are the only pairs with the same number of each
    of the above types of edges.

    Since the right edge from the base of the graph is solid, our
    candidate for $Q_{1}$ must have a solid right edge.  Observe that
    all paths in the left half of this graph with at least one solid
    right edge have two dashed right edges.  The only candidates for
    $Q_{2}$ with two dashed edges are $P_{2}$ and $P'_{2}$, both of
    which have two dashed left edges.  The only remaining candidate
    for $Q_{1}$ with two dashed left edges is in fact $P_{1}$.

\end{proof}


\begin{lemma}
    \label{lemma:deg7}
    Suppose $T$ is a tableau such that $A(T) = B(T)$, $C(T) = D(T)$,
    $e_{1}T\neq 0$, and $e_{2}T\neq 0$.  Then $T$ has a degree 7
    relation above it.
\end{lemma}

\begin{figure}[tbp]
    \centerline{\includegraphics[height=7in]{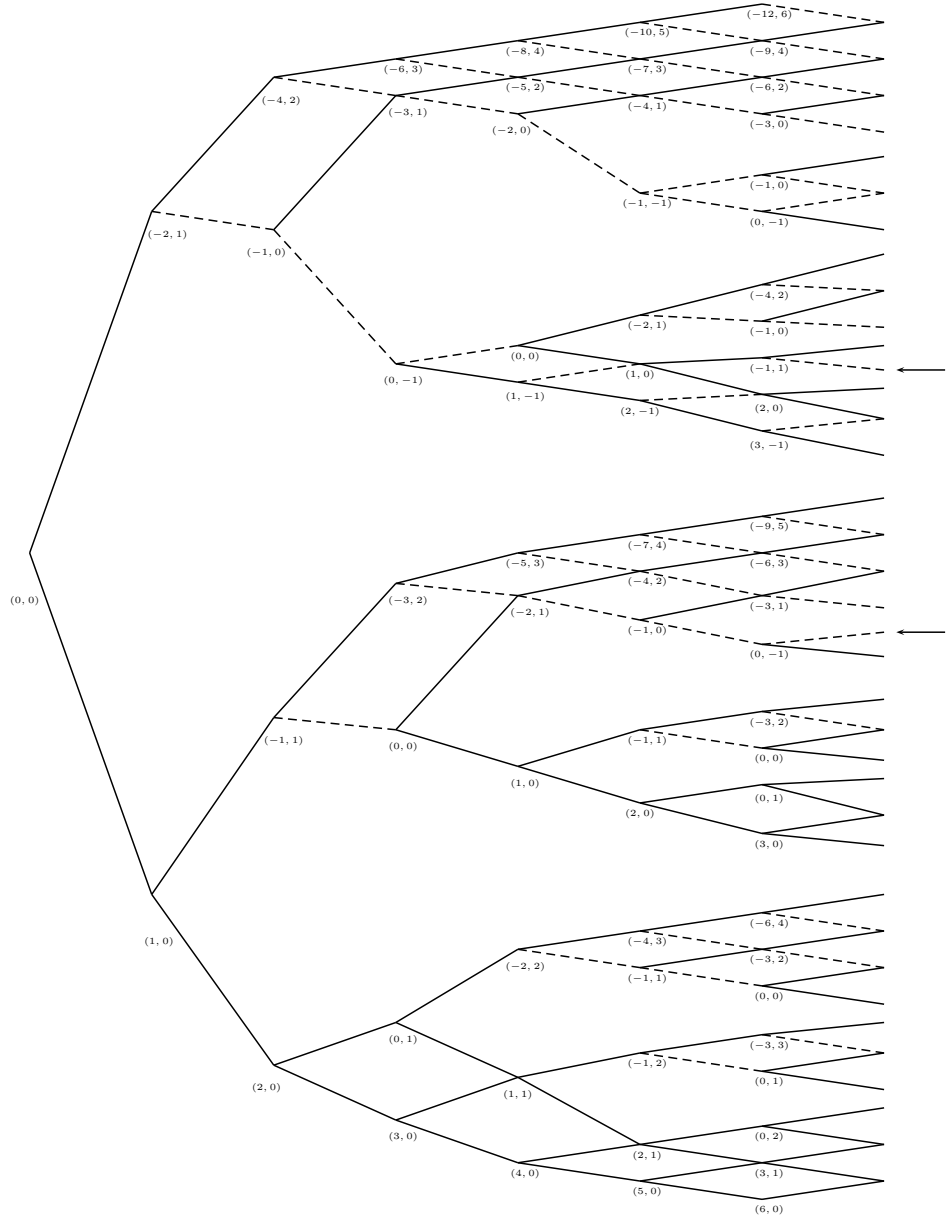}}
    \caption{Picture for Lemma \ref{lemma:deg7}}
    \protect\label{piclem7}
\end{figure}

\begin{proof}

    Note that in order to increase the readability of the graph in
    Figure \ref{piclem7}, it has been oriented to grow to the right
    rather than up.  We therefore take a down edge to indicate acting
    by $e_{1}$ and an up edge to indicate acting by $e_{2}$.
    Otherwise, the legend from Table \ref{piclegend} applies.

    We must first confirm that the sequences
    $e_{2}e_{1}^{2}e_{2}^{3}e_{1}$ and $e_{1}e_{2}^{3}e_{1}^{2}e_{2}$
    do not produce $0$ when applied to $T$.  By Sublemma
    \ref{sublemma:dashed}, we need only show that $e_{1}^{3}e_{2}T$,
    $e_{2}^{2}e_{1}^{3}e_{2}T$, and $e_{2}^{2}e_{1}T$ are not $0$.
    First note that there is at least one $\bar{1}$ in the top row of
    $T$, and the application of $e_{1}^{2}$ to $e_{2}T$ acts on
    entries corresponding to the left block of $+$'s.  It follows that
    the $\bar{1}$'s in the top row of $T$ are also present in
    $e_{1}^{2}e_{2}T$, so $e_{1}^{3}e_{2}T \neq 0$.  This final
    application of $e_{1}$ changes a $\bar{1}$ to a $\bar{2}$.  Since
    $e_{2}$ acts on the left block of $+$'s in $e_{1}^{3}e_{2}T$, it
    leaves this $\bar{2}$ alone, and it can be acted on by the next
    application of $e_{2}$, so $e_{2}^{2}e_{1}^{3}e_{2}T \neq 0$.
    Finally, note that there is a $\bar{2}$ in the top row of $T$ and
    $e_{1}$ changes a $\bar{1}$ to a $\bar{2}$ in the top row of $T$.
    Thus, there are at least two $\bar{2}$'s in the top row of
    $e_{1}T$, and $e_{2}^{2}e_{1}T \neq 0$.

    Now that we know that neither of these sequences produces $0$ when
    applied to $T$, it is clear that we have
    $e_{2}e_{1}^{2}e_{2}^{3}e_{1}T = e_{1}e_{2}^{3}e_{1}^{2}e_{2}T$,
    since the paths corresponding to these sequences leading from the
    base of the graph in Figure \ref{piclem7} to the leaves marked by
    arrows have one solid down edge, three dashed down edges, two
    solid up edges, and one dashed up edge.  Note that these paths are
    equivalent to $e_{2}e_{1}e_{2}e_{1}e_{2}^{2}e_{1}T =
    e_{1}e_{2}^{2}e_{1}e_{2}e_{1}e_{2}T$, due to the degree 2
    relations above $e_{2}^{2}e_{1}T$ and $e_{1}e_{2}T$; we denote
    these alternative paths by $P'_{1}$ and $P'_{2}$ respectively.  We
    must now confirm that among all pairs $(Q_{1},Q_{2})$ of
    increasing paths from the base in these graphs such that $Q_{1}$
    begins by following the up edge and $Q_{2}$ begins by following
    the down edge, $(P_{1}, P_{2})$, $(P_{1},P'_{2})$,
    $(P'_{1},P_{2})$ and $(P'_{1},P'_{2})$, are the only pairs with
    the same number of each of the above types of edges.

    We first address pairs of paths of length no greater than 5.  For
    a path to be a candidate for $Q_{1}$, it must have at least one
    solid down edge.  The only such paths are those beginning with
    $e_{1}^{2}e_{2}T$.  As these paths have two dashed down edges,
    their only possible $Q_{2}$ mate is $e_{1}^{2}e_{2}^{2}e_{1}T$,
    but none of our $Q_{1}$ candidates have the same edge content as
    this path.
    
    We now consider paths of length 6.  As in the preceding paragraph,
    our only candidates for $Q_{1}$ are those paths that contain a
    solid down edge and begin with $e_{1}^{2}e_{2}T$; all such paths
    have exactly two dashed down edges.  Up to degree 2 relations,
    there are three candidates for $Q_{2}$:
    $e_{1}^{2}e_{2}^{3}e_{1}T$, $e_{1}e_{2}e_{1}^{2}e_{2}e_{1}T$, and 
    $e_{1}^{2}e_{2}^{2}e_{1}^{2}T$.  None of these paths contain a
    dashed up edge, which leaves only $e_{1}^{5}e_{2}T$ as our only
    candidate for $Q_{1}$; this cannot be paired with any of our three
    potential $Q_{2}$ paths.

    Finally, we restrict our attention to paths of length 7.  There
    are six paths (again, up to degree 2 relations) in the top half of
    the graph with solid down edges: $e_{1}^{5}e_{2}^{2}T$ and those
    paths beginning with $e_{1}^{3}e_{2}T$.  The former has four
    dashed down edges, a feature lacking from all paths in the bottom 
    half of the graph.  We may also exclude from our consideration
    $e_{1}^{6}e_{2}T$, as all candidates for $Q_{2}$ with only one up 
    edge have at most one dashed down edge.
    
    The remaining three paths that might be $Q_{1}$ all have a dashed 
    up edge; the only $Q_{2}$ candidates with this feature are $P_{2}$
    and $P'_{2}$.  The only paths in the top half of the graph with
    the same edge content as these are $P_{1}$ and $P'_{1}$.

\end{proof}

\begin{example}
    
    In Figure \ref{fig:lem7}, we have a crystal in which the bottom
    tableau $T$ has the statistics $A(T) = B(T) = 1$ and $C(T) =
    D(T) = 0$, illustrating Lemma \ref{lemma:deg7}.
    
    \begin{figure}[!t]
        \centerline{\includegraphics[height=4.55in]{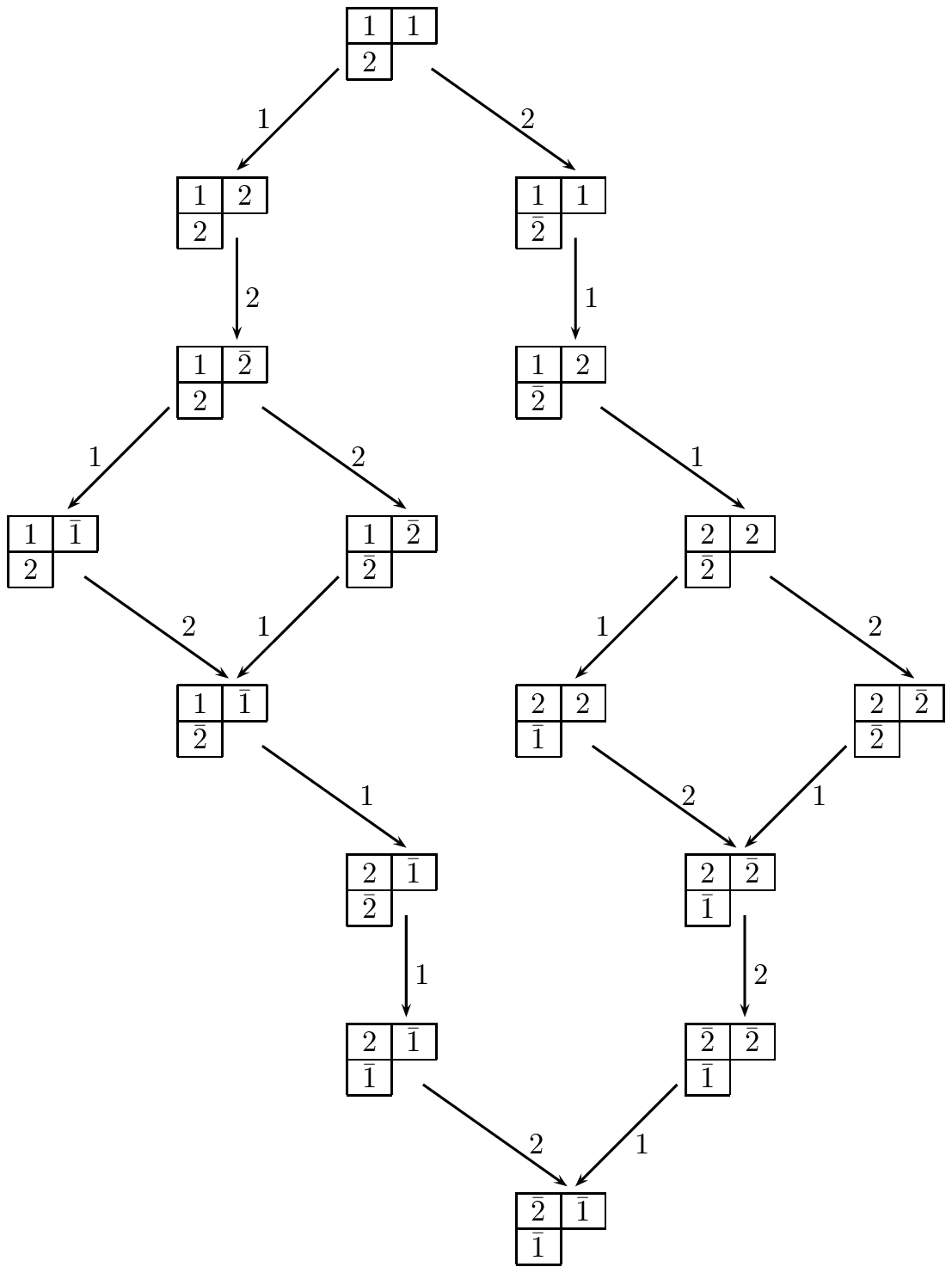}}
        \caption{}
        \protect\label{fig:lem7}
    \end{figure}
\end{example}

\section{Further work}

In the program to locally characterize crystal graphs, two questions
immediately arise following this result.  First, can a local
characterization be provided for doubly laced crystals?  And second,
could such a result be provided for triply laced crystals (i.e., those
of type $G_{2}$)?

It is very reasonable to suspect that a set of local graph theoretic
axioms that characterize doubly laced crystals exists.  It appears
that they may need to be ``less local'' than the axioms in \cite{Stem}
for simply laced crystals.  For instance, we have seen that when $T$
has a degree 5 relation above it, there is a degree 2 relation above
$e_{2}e_{1}T$.  Thus, one of these axioms might be of the form ``If
$v$ is a vertex satisfying certain local conditions, then $v$ must
have a degree 5 relation above it and $e_{2}e_{1}v$ must have a degree
2 relation above it.

It may be possible to prove that such a set of axioms characterize
doubly laced crystals by using virtual crystals \cite{OSS}, a
construction that realizes non-simply laced crystals in terms of
embeddings into simply laced crystals.  More precisely, one can
construct a ``virtualization'' of each of the relations dealt with
above; each of these would be a local piece of a type $A_{3}$ crystal
that corresponds to these relations in terms of the virtual crystal
embeddings.  It would then suffice to show that when these virtual
pieces are assembled according to the doubly-laced axioms, the simply
laced axioms are satisfied.

Calculations suggest that there are over 40 different relations in the
case of $G_{2}$ crystals, some of degree greater than 10
\cite{StemPC}.  The methods employed here would clearly be inadequate
to produce a human-readable proof of a local description of such
graphs.  However, there are probably statistics on $G_{2}$ Littelmann
paths similar to the $ABCD$ statistics used here that could be used to
reduce the problem to a finite number of cases; these could, in turn,
be checked by computer.

\end{document}